\documentclass[conf]{new-aiaa}
\usepackage[utf8]{inputenc}

\usepackage{graphicx}
\usepackage{amsmath}
\usepackage[version=4]{mhchem}
\usepackage{siunitx}
\usepackage{longtable,tabularx}
\usepackage{todonotes}
\usepackage{caption}
\usepackage{subcaption}
\usepackage{multirow}
\title{Heuristic for Min-Max Heterogeneous Multi-Vehicle Multi-Depot Traveling Salesman Problem}

\author{Deepak Prakash Kumar\footnote{Graduate Student, Mechanical Engineering, 3123 TAMU, College Station, TX 77843.}, Sivakumar Rathinam\footnote{Professor, Mechanical Engineering, 3123 TAMU, College Station, TX 77843.}, Swaroop Darbha\footnote{Professor, Mechanical Engineering, 3123 TAMU, College Station, TX 77843.}}
\affil{Texas A\&M University}
\author{Trevor Bihl\footnote{Senior Research Engineer, Sensing Management Branch, Dayton, OH 45433.}}
\affil{Air Force Research Laboratory}

\begin{document}

\maketitle

\begin{abstract}
In this article, a heuristic is proposed for a min-max heterogeneous multi-vehicle multi-depot traveling salesman problem (TSP), wherein heterogeneous vehicles start from given depot positions and need to cover a given set of targets. The vehicles should cover given targets such that the maximum tour time is minimized. In the considered problem, vehicles considered can be functionally heterogeneous, wherein specific targets must be covered by a particular vehicle, structurally heterogeneous, wherein the vehicles can travel at different speeds, or both. The proposed heuristic generalizes the MD heuristic for the min-max homogeneous multi-vehicle multi-depot TSP and has three stages: an initialization stage to generate a feasible solution, a local search stage in which the vehicle with the maximum tour time is improved, and a perturbation stage to break from a local minimum. The proposed heuristic is benchmarked with the optimal solution obtained by solving a mixed integer linear program using branch and cut for instances considering three vehicles covering thirty targets. Variations in the percentage of vehicle-target assignments and the number of vehicles starting at the same depot are studied to show the heuristic's effectiveness in producing high-quality solutions. It was observed that the heuristic generated feasible solutions within $4\%$ of the optimum on average for the considered instances.
\end{abstract}

\section{Nomenclature}

{\renewcommand\arraystretch{1.0}
\noindent\begin{longtable*}{@{}l @{\quad=\quad} l@{}}
TSP & Traveling Salesman Problem \\
% MMHTSP & Min-Max Heterogeneous multi-vehicle multi-depot Traveling Salesman Problem \\
MILP & Mixed Integer Linear Program \\
LP & Linear Program
\end{longtable*}}

\section{Introduction}

In vehicle routing, one of the prevalent problems studied in the literature is the traveling salesman problem (TSP), wherein a salesman (or vehicle) needs to cover a set of targets in a graph by visiting each target exactly once at a minimum tour length and return to the starting location \cite{TSP_book_applegate}. An extension of such a problem is a multi-vehicle TSP, wherein multiple vehicles/salesmen must cover given targets starting from and ending at given depots. Two variations of such a problem have been studied in the literature that differ based on the objective function considered: the min-sum and the min-max variants. While the min-sum variant aims to construct tours such that the total cost of the vehicles is minimized, the min-max variant minimizes the maximum tour cost. One of many applications for studying such problems includes surveillance applications utilizing unmanned aerial or ground vehicles.
% which could start at different locations (depots). 
For such applications, it is essential to construct tours for functionally and structurally heterogeneous vehicles that could start at different depots. It should be noted that functional heterogeneity of vehicles can arise due to different sensing capabilities. On the other hand, structural heterogeneity could arise from vehicles traveling at different speeds or having different motion constraints, such as in the case of a Dubins vehicle \cite{Dubins}.

Efficient algorithms for the min-sum heterogeneous multi-vehicle multi-depot TSP have been developed in the literature. For example, a graph transformation is utilized in \cite{todays_TSP} to transform the problem to a single-vehicle asymmetric TSP and utilize well-known single-vehicle TSP solvers, such as LKH \cite{LKH_TSP_solver}. Alternatively, in \cite{algorithms_heterogeneous_multi_depot_multi_vehicle}, the min-sum variant was shown to be solved to optimality for instances with 100 targets and five vehicles using branch and cut in 300 seconds on average. 
% In \cite{algorithms_heterogeneous_multi_depot_multi_vehicle}, a branch and cut algorithm is used to plan tours for vehicles that are functionally heterogeneous and with different motion constraints, and instances with 100 targets and 5 vehicles were solved in 300 seconds due to a tight lower bound from the LP solution.
A tight lower bound is obtained from the relaxed linear program (LP) solution for the min-sum variant, as is demonstrated in \cite{algorithms_heterogeneous_multi_depot_multi_vehicle}. However, for the min-max variant, the lower bound obtained is weak in comparison. Hence, the min-max variant is typically addressed using a heuristic or an approximation algorithm.

The min-max variant of the problem has been studied under two variations: the single-depot variant and the multi-depot variant. While many studies have considered heuristics for the single-depot variant, a review of which is discussed in \cite{survey_mtsp, Memetic_search_min_max_TSP}, fewer studies have focused on the multi-depot variant of the problem. The min-max multi-vehicle multi-depot TSP for homogeneous vehicles was first proposed in \cite{min_max_vrp_LB_based_load_balancing}, and the authors proposed four heuristics for the same. They observed that two heuristics based on an LP-based load balancing heuristic and a region partition heuristic performed the best for 13 instances with uniform customer locations. An ant colony-based optimization was proposed in \cite{ant_colony_min_max} and was compared with the LP-based heuristic for a few instances for targets obtained from a uniform distribution. The proposed heuristic performed better than the LP-based heuristic. In \cite{MD_algorithm}, the authors proposed a heuristic called the "MD" algorithm, wherein an initial solution was first obtained from LP-based load balancing in \cite{min_max_vrp_LB_based_load_balancing}. The obtained solution was improved using a local search to improve the maximal route by removing a target from the maximal route and inserting it in another vehicle's route. Finally, a perturbation stage was used to break from a local minimum. The authors benchmarked the proposed heuristic over 43 instances against three other heuristics considering varying target-to-vehicle ratios, varying number of vehicles at a depot, and different distributions for generating the targets. An improvement to the MD algorithm was studied in \cite{min_max_split_delivery_vrp} as a part of a min-max split delivery problem, wherein additional neighborhoods, such as cyclic transfer of targets and two-point swap, were included. Recently, a memetic algorithm was proposed in \cite{Memetic_search_min_max_TSP}, wherein multiple initial solutions were generated, using which a variable neighborhood descent was performed over six types of neighborhoods to provide the best-known solutions for many of the instances considered in \cite{MD_algorithm}.
% In the first heuristic, an LP was solved for initial assignment, the Concorde TSP solved was used to generate tours, and tours were improved by changing the number of customer assignments to each vehicle., wherein the region was partitioned into equal convex polygons 
% \cite{min_max_tsp_structural_heterogeneous_Dubins}

While the above-discussed studies focus on the min-max routing for homogeneous vehicles, only some studies have explored the problem with multiple heterogeneous vehicles. In \cite{3_approx_algo_2_vehicle_TSP}, a $\frac{3k}{2}$ approximation ratio algorithm was proposed, where $k$ represents the number of vehicles. A feasible solution based on rounding the LP solution was constructed and shown to be at most $\frac{3k}{2}$ times the LP solution. In \cite{heuristic_struct_het_distict_dep}, the authors develop a heuristic based on a primal-dual algorithm for a min-max TSP for vehicles with different average speeds. However, the authors assume that the vehicles must start at distinct depots, and all vehicles must have distinct speeds.
%
% Further, it is also important It should be noted that functional heterogeneity of vehicles can arise for vehicle-target compatibility constraints, and structural heterogeneity of vehicles could arise from vehicles traveling at different speeds or having different motion constraints, such as in the case of a Dubins vehicle \cite{Dubins}.

From the discussed studies, it can be observed that while an approximation algorithm is presented in \cite{3_approx_algo_2_vehicle_TSP} for a min-max heterogeneous multi-vehicle multi-depot TSP considering functional and structural heterogeneity, a heuristic for providing a good feasible solution for this problem has not been studied. Hence, in this paper, a heuristic is proposed for this problem by considering vehicles to be structurally heterogeneous due to different vehicle speeds in addition to functional heterogeneity. Further, the presented heuristic is implemented and compared with the optimal solution obtained from solving the mixed-integer linear program (MILP) for this problem using branch and cut for instances with three vehicles and thirty targets. Further, the fraction of targets pre-assigned to vehicles, and the number of vehicles starting at a depot are varied to study the quality of the feasible solution obtained from the heuristic. It was observed that the feasible solution deviated from the optimal solution by about $4\%$ on average for the instances considered.

\section{Problem Description}

In this study, $k$ vehicles are considered with vehicle speeds $v_1, \cdots, v_k$ that start at depots $d_1, d_2, \cdots d_k.$ It should be noted that the considered depots for vehicles and the vehicle speeds need not necessarily be distinct. Let $T$ denote all targets that must be covered by the $k$ vehicles. Further, let $R_i$ denote the subset of targets in $T$ that must be necessarily covered by the $i\textsuperscript{th}$ vehicle. It should be noted that $R_i \cap R_j = \emptyset \forall i, j \in \{1, 2, \cdots, k\}$ for $i \neq j$, and $\bigcup_{i} R_i \subseteq T$. It is desired to construct tours for the $k$ vehicles such that 
\begin{itemize}
    \item Vehicle $i$ covers all targets in $R_i$ exactly once,
    \item Each target in $T \setminus \bigcup_{i} R_i$ is covered by exactly one vehicle once, and
    \item The maximum tour time of the $k$ vehicles is minimized.
\end{itemize}
Let $E_i$ denote the set of all edges that connect vertices in $T \bigcup \{d_i\}$ for the $i\textsuperscript{th}$ vehicle. Assuming that the vehicles travel at a constant speed, the time taken to travel from vertex $m$ to $n$ for the $i\textsuperscript{th}$ vehicle is equal to the Euclidean distance between $m$ and $n$ divided by the vehicle speed $v_i$ for all $(m, n) \in E_i$. Hence, the graph $G_i = (T \bigcup \{d_i\}, E_i)$ corresponding to the $i\textsuperscript{th}$ vehicle is complete, undirected, symmetric, and satisfies the triangle inequality. The MILP formulation for the considered problem is given in \cite{3_approx_algo_2_vehicle_TSP}.

\section{Heuristic Description}

The proposed heuristic in this study to address the described routing problem is inspired by the MD algorithm \cite{MD_algorithm}. The MD algorithm was proposed for a min-max homogeneous multi-vehicle multi-depot TSP and utilized three stages: an initialization stage using load balancing, local search to improve the maximal tour, and a perturbation stage to break from a local minimum. The heuristic proposed in this study similarly uses three stages to obtain a feasible solution.

\subsection{Initialization}

In this stage, an initial feasible solution is obtained using a modified load balancing approach used in \cite{MD_algorithm}. Since the vehicles are considered to travel at different speeds, each vehicle is desired to be allocated a fraction of the number of targets depending on its speed. Hence, vehicle $i$ is desired to be allocated about $|T| \frac{v_i}{\sum_{i = 1}^k v_i}$ number of targets from $T$. Since vehicle $i$ needs to necessarily cover targets in $R_i \subset T$, vehicle $i$ is desired to be allocated at least $\Big{\lfloor} |T| \frac{v_i}{\sum_{i = 1}^k v_i} \Big{\rfloor} - |R_i|$ number of targets from the set $T \setminus \bigcup_i R_i,$ i.e., the set of targets that any of the $k$ vehicles can cover. It should be noted that the floor operator was used to ensure that the number of targets to be allocated to the $i\textsuperscript{th}$ vehicle is an integer and that a feasible solution can always be obtained.

An integer program is considered to determine the allocation of targets. Consider a binary variable $x_{tj}$ for $t \in T \setminus \bigcup_i R_i, j = 1, 2, \cdots, k$ which denotes the allocation of target $t$ to vehicle $j$. It should be noted that if $x_{tj} = 1,$ target $t$ is allocated to vehicle $j$. The integer program formulated for the allocation of targets in $T \setminus \bigcup_i R_i$ is as follows:
\begin{align}
    &\min \sum_{t \in T \setminus \bigcup_i R_i} \sum_{j = 1}^n c_{tj} x_{tj}, \\
    \text{s.t.} \quad & \sum_{j = 1}^n x_{tj} = 1, \quad \forall t \in T \setminus \bigcup_i R_i, \label{eq: target_allocation} \\
    & \sum_{t \in T \setminus \bigcup_i R_i} x_{tj} \geq \bigg{\lfloor} |T| \frac{v_j}{\sum_{i = 1}^k v_i} \bigg{\rfloor} - |R_j|, \quad \forall j = 1, 2, \cdots, k, \label{eq: vehicle_num_targets_assigned} \\
    & x_{tj} \in \{0, 1 \} \quad \forall t \in T \setminus \bigcup_i R_i, \,\, j = 1, 2, \cdots, k,
\end{align}
where $c_{tj}$ represents the time taken to travel from depot $d_j$ to target $t$ for vehicle $j$. In the above formulation, constraint~\eqref{eq: target_allocation} ensures that each target in $T \setminus \bigcup_i R_i$ is allocated to exactly one vehicle, and constraint~\eqref{eq: vehicle_num_targets_assigned} specifies the minimum number of targets allocated to vehicle $j$. An allocation of targets in $T \setminus \bigcup_i R_i$ is obtained by solving the relaxation of the above program and rounding. A tour for the $i\textsuperscript{th}$ vehicle is then constructed using LKH \cite{LKH_TSP_solver} to cover the obtained allocation of targets in $T \setminus \bigcup_i R_i$ and targets in $R_i$ for all $i = 1, \cdots, k$.

It should be noted that when multiple vehicles are present at the same depot, the initial location of each vehicle that starts at the same location is perturbed, similar to \cite{MD_algorithm} to obtain an initial feasible solution from the load balancing. The perturbation accounts for when multiple vehicles with the same speed start from the same initial location. In this particular case, varying target allocations among the multiple vehicles that start at the same location does not change the objective function value without the perturbation. For the implementation, the initial locations were placed symmetrically around a circle of radius $0.1$ centered at the initial location (before perturbation). The angle of perturbation for the first vehicle was chosen randomly.

\subsection{Local search}

A local search is performed to improve the tour of the maximal vehicle, i.e.,  the vehicle with the maximum tour time. For this purpose, a target is removed from the vehicle's tour and is inserted into another vehicle's tour, similar to \cite{MD_algorithm}. A metric called the savings is defined for each target that could be removed from the maximal tour. The metric estimates the reduction in tour time for the maximal vehicle by removing the corresponding target from the tour. Denoting the vehicle with the maximum tour time with index $i$, the savings associated with removing target $t$ is defined as
\begin{align}
    \text{savings}_{t, i} = \frac{distance(t_{prev}, t)}{v_i} + \frac{distance(t, t_{next})}{v_i} - \frac{distance(t_{prev}, t_{next})}{v_i},
\end{align}
where $t_{prev}$ and $t_{next}$ are vertices covered before and after target $t$ by vehicle $i$. The above-defined savings is computed for all targets covered by vehicle $i$ except for targets in the set $R_i$ since targets in $R_i$ should be necessarily covered by vehicle $i$. The targets covered by vehicle $i$ (that are not in $R_i$) are ranked in descending order of their savings. 
% It should be noted here that the computed savings are conservative. This is because the obtained tour for vehicle $i$ by removing $t$ and connecting $t_{prev}$ and $t_{next}$ is a feasible tour for vehicle $i$, as shown in Fig.~ .
% The representation of the corresponding tour obtained using the proposed removal of $t$ is shown in Fig.~ .

It is desired to remove the target with the highest savings first. For this purpose, the vehicle to which it must be inserted is identified by defining a ``cost increase" metric. The metric estimates the increase in tour time for the considered vehicle by inserting the corresponding target in the considered vehicle's tour. Suppose vehicle $j$'s current tour is given by $(u_0 = d_j, u_1^j, u_2^j, \cdots, u_{n_j - 1}^j, u_{n_j}^j = d_j)$. The cost associated with inserting a target $t$ between vertices $u_l$ and $u_{l + 1}$ is defined as
\begin{align}
    \text{cost increase}_{t, j, (u_l^j, u_{l + 1}^j)} = \frac{distance(u_l^j, t)}{v_j} + \frac{distance(t, u_{l + 1}^j)}{v_j} - \frac{distance(u_l^j, u_{l + 1}^j)}{v_j}.
\end{align}
Since target $t$ can be inserted between any pair of vertices in vehicle $j$'s tour, the cost increase for vehicle $j$ for inserting target $t$ is defined to be the minimum among all such costs. The defined cost increase can be obtained corresponding to target $t$ for all vehicles $j \in \{1, 2, \cdots, k\} \setminus \{i\}.$ The vehicle with the least cost of increase is chosen for inserting target $t$, denoted with index $m$. 
% It should be noted that as opposed to ``savings", ``cost increase" for each vehicle is an upper bound for the increase in cost of the corresponding vehicle's tour. This is because a feasible solution is obtained by inserting $t$ in 
% The representation of the corresponding tour obtained using the proposed insertion of $t$ is shown in Fig.~ .

Finally, LKH is used to obtain tours for vehicles $i$ and $m$, wherein target $t$ is removed from vehicle $i$ and allocated to vehicle $m$. If the new maximal tour's cost is less than the previous one, the solution is updated. If not, the target with the second highest savings previously identified is attempted to be removed from vehicle $i$ using the same process as above. The search is terminated when all targets on the maximal tour have been considered and no improvement in the solution is obtained.

\subsection{Perturbation}

\begin{figure}[htb!]
    \centering
    \includegraphics[width = 0.3\linewidth]{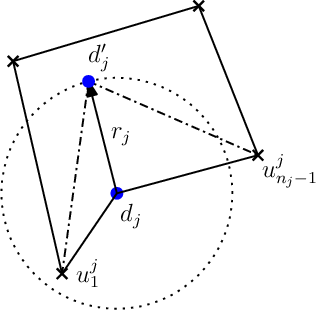}
    \caption{Depiction of perturbation for the $j\textsuperscript{th}$ vehicle}
    \label{fig: perturbation_j_vehicle}
\end{figure}

From the local search, a local minimum is obtained and cannot be improved by removing a target from the maximal vehicle and reassigning it to another vehicle. The obtained solution is perturbed using a similar method in the MD algorithm \cite{MD_algorithm}. For the $j\textsuperscript{th}$ vehicle, the average cost to the targets connected to the vehicle's depot is computed and is given by $r_j = \frac{1}{2} \left(\frac{distance(d_j, u_1^j) + distance(d_j, u_{n_j - 1}^j)}{v_j} \right)$ for $j = 1, \cdots, k$. Here, $u_1^j$ and $u_{n_j - 1}^j$ denote the targets connected to $d_j$.
% . It should be noted that since the cost of traveling between vertices $m$ to $n$ for vehicle $i$ is given by $\frac{distance(m, n)}{v_i}$, the distance between two vertices can always be obtained to compute the required average distance. 
Then, in the first round of perturbation, each vehicle's depot is perturbed from its initial location by an arbitrary angle to a location that is at the obtained average distance, as shown in Fig.~\ref{fig: perturbation_j_vehicle}. Using the same allocation of targets from the solution obtained from the local search, the tours for the vehicles are computed using LKH in the graph with the perturbed depot location to obtain a feasible solution. A local search is then performed on this graph with perturbed depot locations to obtain a new set of targets allocated for each vehicle. The tours for the vehicles are then constructed using LKH in the original graph (with depots at their initial locations) to obtain a feasible solution different from the solution obtained before perturbation. If an improved solution is obtained compared to the solution before perturbation, then the solution is updated.

Similar to \cite{MD_algorithm}, the perturbation step is repeated until there is no improvement for five consecutive iterations. The perturbation angle for a depot in a subsequent iteration was chosen to be at an angle of $144^\circ$ from the previous iteration. This ensures that the angle of perturbation of a depot in the sixth consecutive iteration without improvement is the same as the first iteration.

\section{Results}

The presented heuristic was implemented in Python~3.9 on a laptop with AMD Ryzen $9$ $5900$HS CPU running at $3.30$ GHz with $16$ GB RAM. Further, to benchmark the presented heuristic, the mixed integer linear program formulation given in \cite{3_approx_algo_2_vehicle_TSP} was implemented using branch and cut and solved using Gurobi \cite{gurobi} in Python. To this end, subtour elimination constraints were used as cutting planes and lazy cuts using lazy callback functions. Further, the feasible solution obtained from the presented heuristic was used as an initial incumbent solution for branch and cut. In addition, the solution obtained from rounding the LP solution, which corresponds to a solution with an approximation ratio of $\frac{3k}{2}$, where $k$ represents the number of vehicles, was also obtained. The obtained feasible solution was compared with the initial feasible solution obtained through load balancing.

% In \cite{3_approx_algo_2_vehicle_TSP}, a  approximation algorithm is presented for the considered problem, . In \cite{3_approx_algo_2_vehicle_TSP}, the feasible solution was obtained by rounding the obtained linear program (LP) solution. In this study, the obtained feasible solution is also compared with the initial feasible solution obtained through load balancing.

In this study, three vehicles with thirty targets were considered. The targets were generated randomly on a $200\times200$ grid. Two scenarios were considered, wherein all the vehicles either start at different depots, or two vehicles with the same speed start at the same depot. The second scenario was studied to evaluate the quality of the initial feasible solution obtained from load balancing. Further, in both scenarios, the fraction of targets for vehicle-target assignments was varied between $0\%, 10\%,$ and $20\%$ to observe the solution quality from the heuristic. For instances with vehicle-target assignments, the targets to be pre-assigned were chosen randomly from the available list of targets, and each such target was randomly allocated to one of the three vehicles. In the implementation, the random.random() function was utilized in Python to generate targets and depots and vehicle-target assignments.

\subsection{Scenario 1: Vehicles start at distinct depots}

\begin{figure}[htb!]
     \centering
     \begin{subfigure}[b]{0.46\textwidth}
         \centering
         \includegraphics[width=\textwidth]{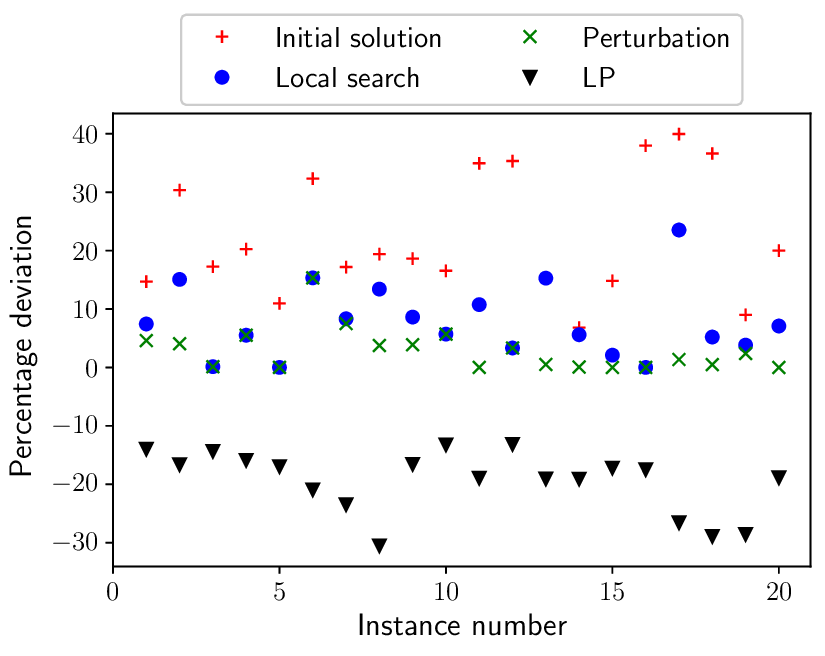}
         \caption{Deviation of heuristic and LP from optimum for $0\%$ target allocation}
         \label{subfig: deviation_heuristic_0_frac_targ_alloc_vehicles}
     \end{subfigure}
     \hfill
     \begin{subfigure}[b]{0.46\textwidth}
         \centering
         \includegraphics[width=\textwidth]{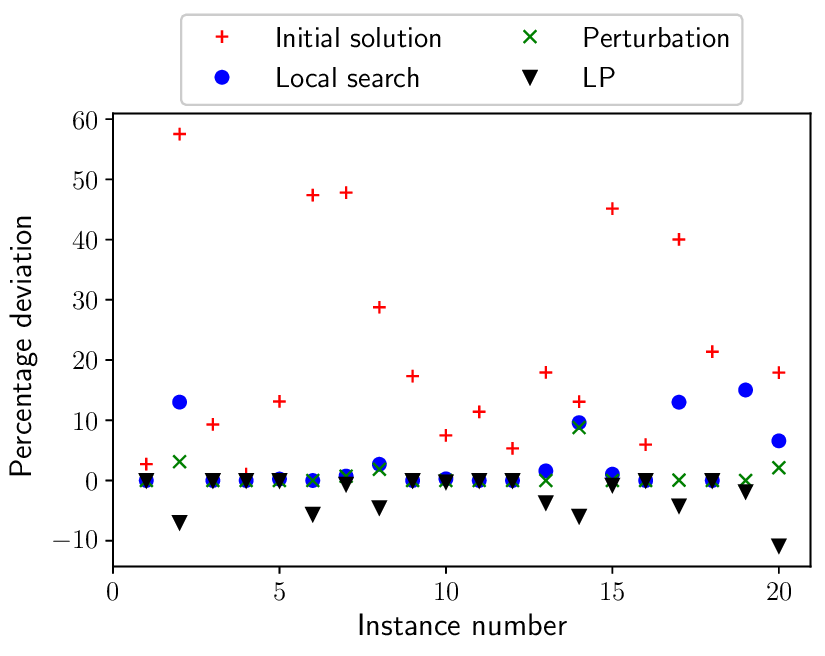}
         \caption{Deviation of heuristic and LP from optimum for $20\%$ target allocation}
         \label{subfig: deviation_heuristic_20_frac_targ_alloc_vehicles}
     \end{subfigure}
     \hfill
     \begin{subfigure}[b]{0.46\textwidth}
         \centering
         \includegraphics[width=\textwidth]{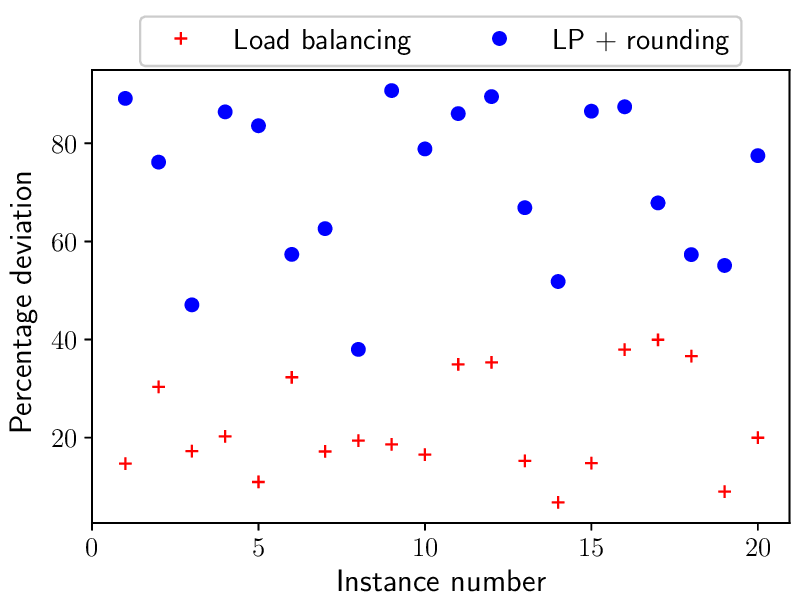}
         \caption{Deviation of two initial feasible solutions from optimum for $0\%$ target allocation}
         \label{subfig: deviation_initial_solution_0_frac_targ_alloc_vehicles}
     \end{subfigure}
     \hfill
     \begin{subfigure}[b]{0.46\textwidth}
         \centering
         \includegraphics[width=\textwidth]{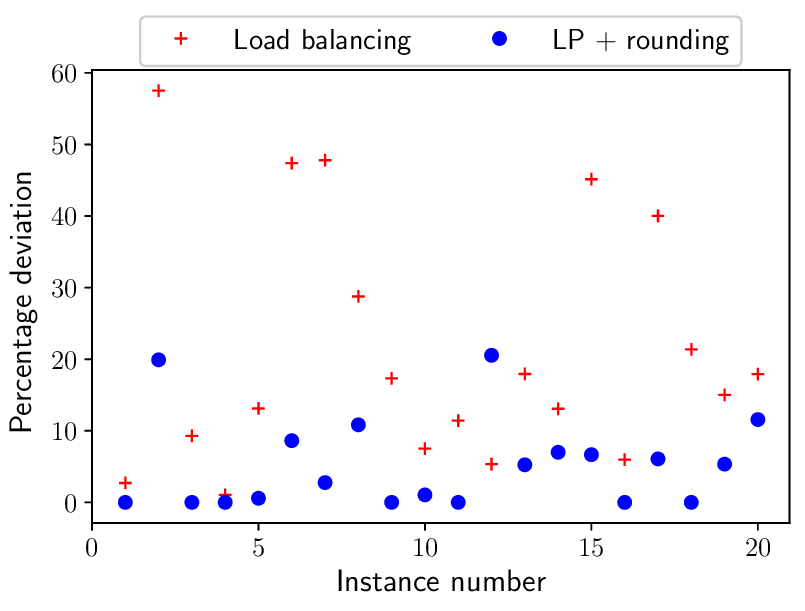}
         \caption{Deviation of two initial feasible solutions from optimum for $20\%$ target allocation}
         \label{subfig: deviation_initial_solution_20_frac_targ_alloc_vehicles}
     \end{subfigure}
     \hfill
     \begin{subfigure}[b]{0.46\textwidth}
         \centering
         \includegraphics[width=\textwidth]{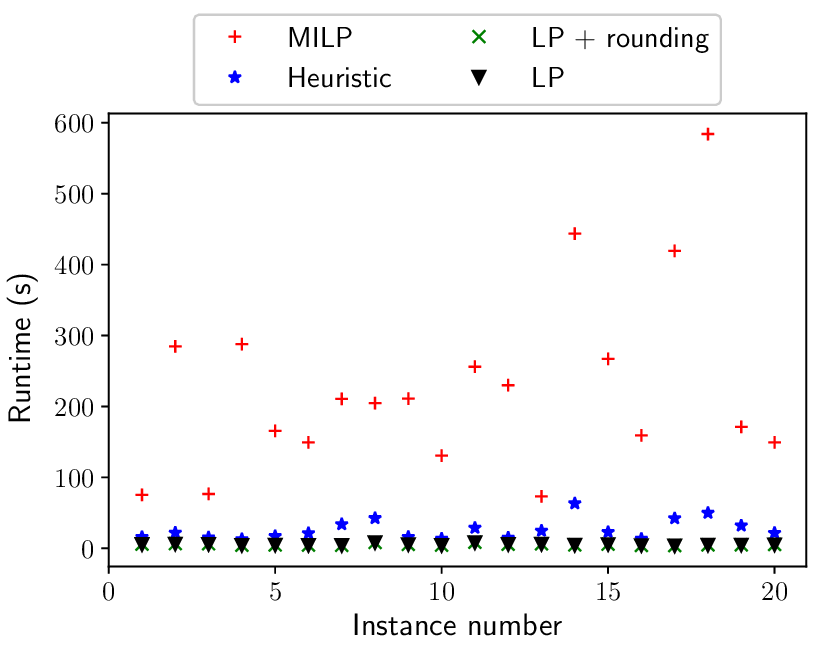}
         \caption{Runtime of branch and cut, heuristic, LP, and LP with rounding for $0\%$ target allocation}
         \label{subfig: runtime_0_frac_targ_alloc_vehicles}
     \end{subfigure}
     \hfill
     \begin{subfigure}[b]{0.46\textwidth}
         \centering
         \includegraphics[width=\textwidth]{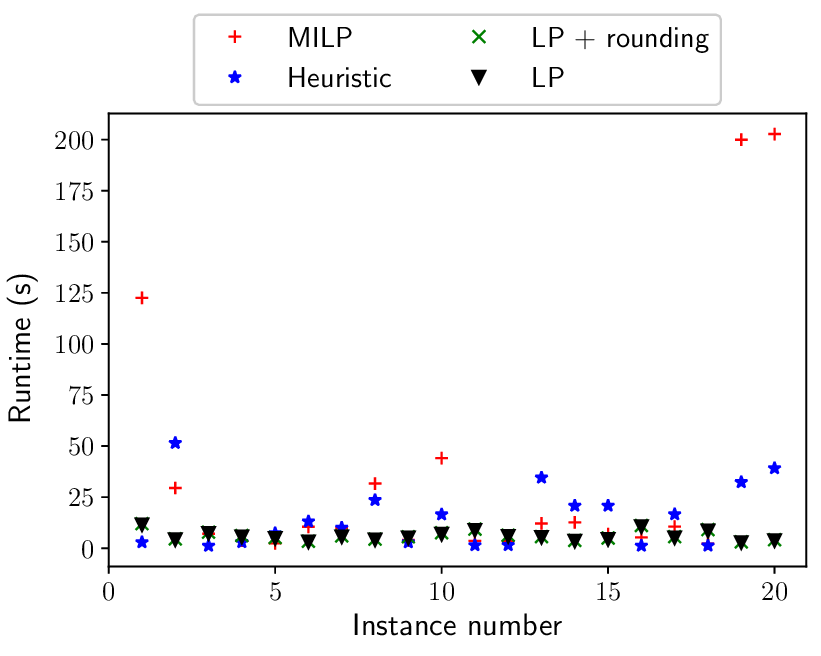}
         \caption{Runtime of branch and cut, heuristic, LP, and LP with rounding for $20\%$ target allocation}
         \label{subfig: runtime_20_frac_targ_alloc_vehicles}
     \end{subfigure}
    \caption{Scenario 1: Vehicles $1, 2,$ and $3$ have speeds $1, 1.5, 2$, respectively, start at different depots covering thirty targets, twenty instances}
    \label{fig: vehicles_diff_depots}
\end{figure}

In this scenario, the speed of the three vehicles was chosen as $1, 1.5,$ and $2$ units, respectively. Further, the initial location of each depot was randomly generated on the same grid. The comparison of the heuristic with the optimal solution obtained using branch and cut for $0\%$ and $20\%$ target assignment and the runtimes are shown in Fig.~\ref{fig: vehicles_diff_depots} for twenty instances. Further, the mean percentage deviation from the optimum for the heuristic, the LP, and the initial feasible solutions are summarized in Table~\ref{tab: mean_percentage_deviation_scenario_1}, and the mean running times are summarized in Table~\ref{tab: mean_running_time_scenario_1}. From Tables~\ref{tab: mean_percentage_deviation_scenario_1} and \ref{tab: mean_running_time_scenario_1} and Fig.~\ref{fig: vehicles_diff_depots}, the following observations can be made:
\begin{itemize}
    \item The final solution obtained from the heuristic was within $4\%$ of the optimal solution on average and was obtained at a lower computation time than solving the MILP. For the scenarios with $0\%, 10\%,$ and $20\%$ target assignment, the number of instances for which the solution obtained from the heuristic was within $2\%$ of the optimal solution was $10, 12,$ and $17$, respectively. It should be noted that the running time for the MILP does not include the computation time for the initial feasible solution.
    \item It can be observed that the lower bound obtained from the LP is poor for the scenario with $0\%$ target assignment but improved with a higher percentage of target assignment. This, in turn, contributed to improved computation time for the MILP for a higher percentage of target assignment.
    \item The initial feasible solution obtained from load balancing performed well for all scenarios. The initial feasible solution from LP with rounding was poor for the scenario with $0\%$ target assignment, since it was observed that many targets were covered by the fastest vehicle in the solution. However, with an increase in target assignment, the solution from the LP with rounding improved significantly.
\end{itemize}

\begin{table}[htb!]
    \centering
    \begin{tabular}{|c|ccccc|c|}
    \hline
    \multirow{3}{*}{\textbf{\begin{tabular}[c]{@{}c@{}}\% target \\ allocation\end{tabular}}} & \multicolumn{5}{c|}{\textbf{Mean percentage deviation from optimum}} & \multirow{3}{*}{\textbf{\begin{tabular}[c]{@{}c@{}}Max. deviation \\ of heuristic from \\ optimum ($\%$)\end{tabular}}} \\ \cline{2-6} 
     & \multicolumn{1}{c|}{LP with} & \multicolumn{1}{c|}{\multirow{2}{*}{Load balancing}} & \multicolumn{1}{c|}{\multirow{2}{*}{Local search}} & \multicolumn{1}{c|}{\multirow{2}{*}{Perturbation}} & \multirow{2}{*}{LP} &  \\
     & \multicolumn{1}{c|}{rounding} & \multicolumn{1}{c|}{} & \multicolumn{1}{c|}{} & \multicolumn{1}{c|}{} & &  \\ \hline
    $0$ & \multicolumn{1}{c|}{$71.80$} & \multicolumn{1}{c|}{$22.41$} & \multicolumn{1}{c|}{$7.81$} & \multicolumn{1}{c|}{$2.93$} & \multicolumn{1}{c|}{$-19.54$} & $15.33$ \\ \hline
    $10$ & \multicolumn{1}{c|}{$16.67$} & \multicolumn{1}{c|}{$29.08$} & \multicolumn{1}{c|}{$7.85$} & \multicolumn{1}{c|}{$3.80$} & \multicolumn{1}{c|}{$-4.28$} & $20.93$ \\ \hline
    $20$ & \multicolumn{1}{c|}{$5.31$} & \multicolumn{1}{c|}{$21.29$} & \multicolumn{1}{c|}{$3.19$} & \multicolumn{1}{c|}{$0.83$} & \multicolumn{1}{c|}{$-2.26$} & $8.79$ \\ \hline
    \end{tabular}
    \caption{Mean percentage deviation of stages of heuristic, LP, and LP with rounding from optimum for $v_1 = 1,$ $v_2 = 1.5$, $v_3 = 2$, thirty targets, and twenty instances}
    \label{tab: mean_percentage_deviation_scenario_1}
\end{table}

\begin{table}[htb!]
    \centering
    \begin{tabular}{|c|c|c|c|c|}
    \hline
    \textbf{\% target allocation} & \textbf{Branch and cut} (s) & \textbf{Heuristic} (s) & \textbf{LP} (s) & \textbf{LP with rounding} (s) \\ \hline
    $0$ & $227.48$ & $26.14$ & $4.97$ & $5.42$ \\ \hline
    $10$ & $93.82$ & $17.84$ & $5.61$ & $6.10$ \\ \hline
    $20$ & $36.23$ & $15.06$ & $5.78$ & $6.31$ \\ \hline
    \end{tabular}
    \caption{Mean running time of branch and cut, heuristic, LP, and LP with rounding for $v_1 = 1,$ $v_2 = 1.5$, $v_3 = 2$, thirty targets, and twenty instances}
    \label{tab: mean_running_time_scenario_1}
\end{table}

\begin{figure}[htb!]
     \centering
     \begin{subfigure}[b]{0.48\textwidth}
         \centering
         \includegraphics[width=\textwidth]{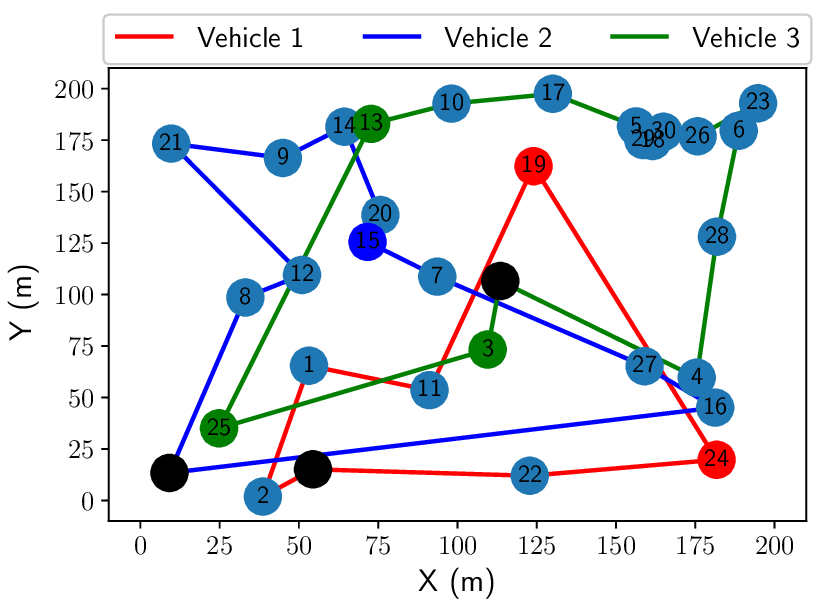}
         \caption{Initial tours obtained from load balancing}
         \label{subfig: initial_tours_load_balancing}
     \end{subfigure}
     % \hfill
     % \begin{subfigure}[b]{0.48\textwidth}
     %     \centering
     %     \includegraphics[width=\textwidth]{Figures/Instance_9_initial_solution_LP_rounding.eps}
     %     \caption{Initial tours obtained from LP with rounding}
     %     \label{subfig: initial_tours_LP_rounding}
     % \end{subfigure}
     \hfill
     \begin{subfigure}[b]{0.48\textwidth}
         \centering
         \includegraphics[width=\textwidth]{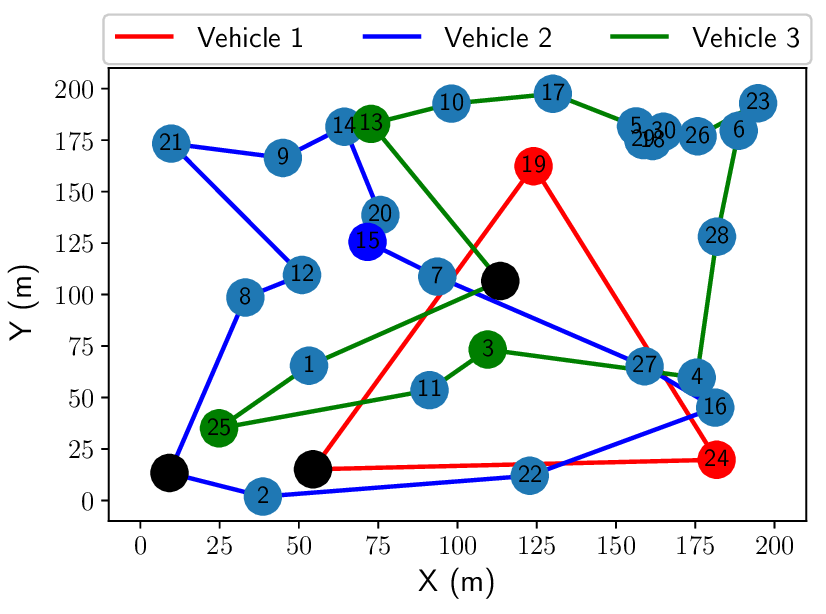}
         \caption{Tours obtained after local search}
         \label{subfig: tours_after_local_search}
     \end{subfigure}
     \hfill
     \begin{subfigure}[b]{0.48\textwidth}
         \centering
         \includegraphics[width=\textwidth]{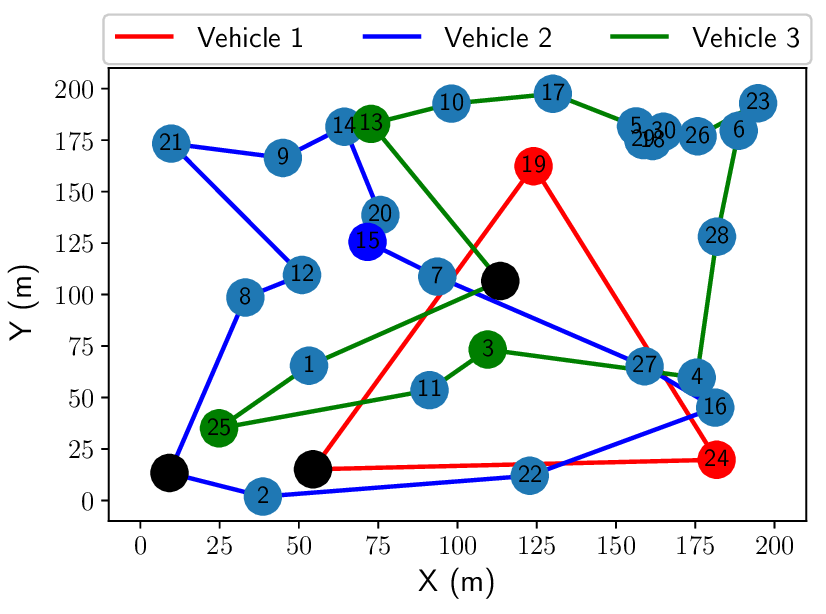}
         \caption{Tours obtained after perturbation}
         \label{subfig: tours_after_perturbation}
     \end{subfigure}
     \hfill
     \begin{subfigure}[b]{0.48\textwidth}
         \centering
         \includegraphics[width=\textwidth]{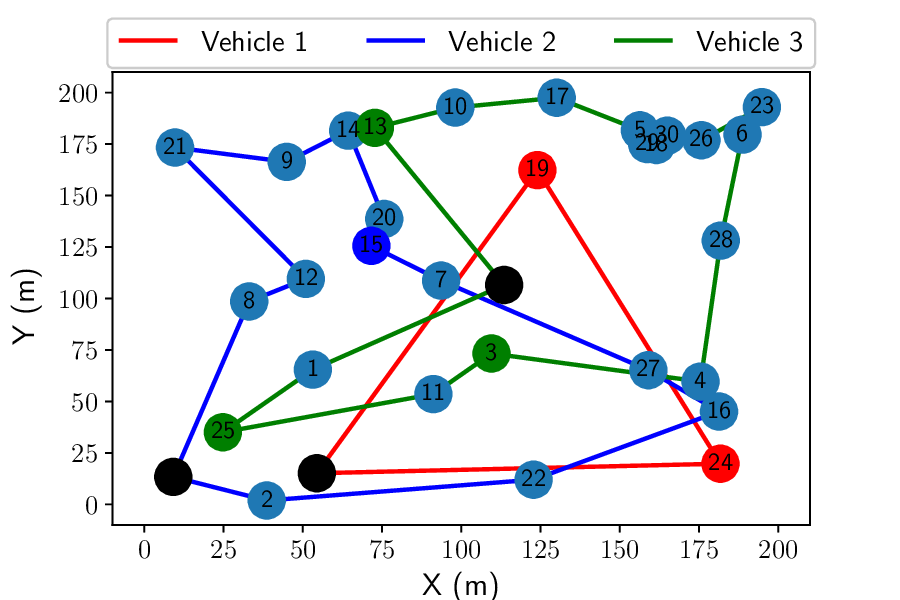}
         \caption{Tours obtained from branch and cut}
         \label{subfig: tours_MILP}
     \end{subfigure}
    \caption{Tours of vehicles for instance 9 corresponding to $v_1 = 1, v_2 = 1.5, v_3 = 2,$ $20\%$ target allocation, and vehicles start at different depots}
    \label{fig: instance_9_tours}
\end{figure}

The tours obtained for the vehicles through the stages of the heuristic and the final solution obtained from solving the MILP for instance $9$ for $20\%$ target allocation are shown in Fig.~\ref{fig: instance_9_tours}. In this figure, the nodes shaded in black represent the depot of each vehicle. Further, targets that must be necessarily covered by a vehicle are shaded in the corresponding color. For this instance, it can be observed that the heuristic obtains the optimal solution after the local search.

\subsection{Scenario 2: Two homogeneous vehicles start at the same depot}

\begin{figure}[htb!]
     \centering
     \begin{subfigure}[b]{0.46\textwidth}
         \centering
         \includegraphics[width=\textwidth]{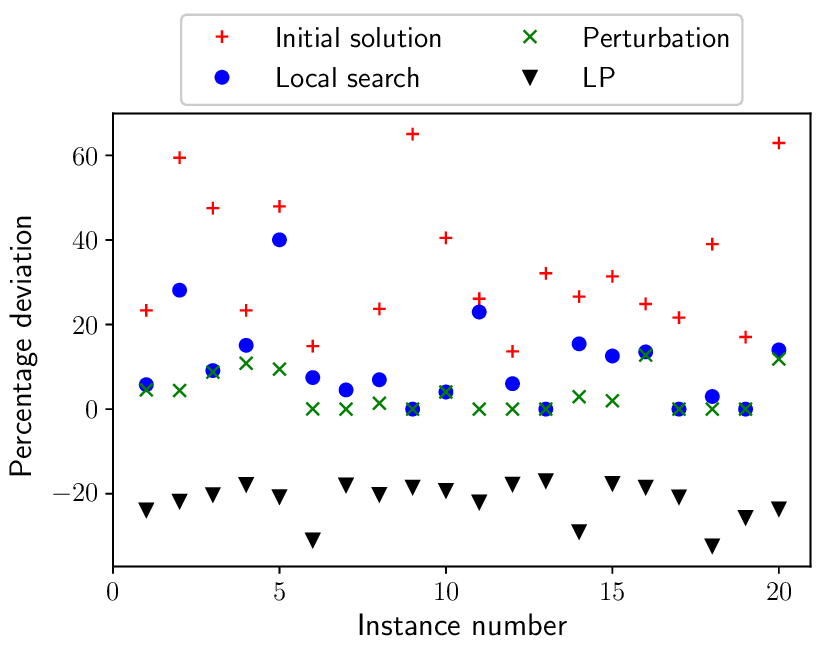}
         \caption{Deviation of heuristic and LP from optimum for $0\%$ target allocation}
         \label{subfig: deviation_heuristic_0_frac_targ_alloc_vehicles_same_depot}
     \end{subfigure}
     \hfill
     \begin{subfigure}[b]{0.46\textwidth}
         \centering
         \includegraphics[width=\textwidth]{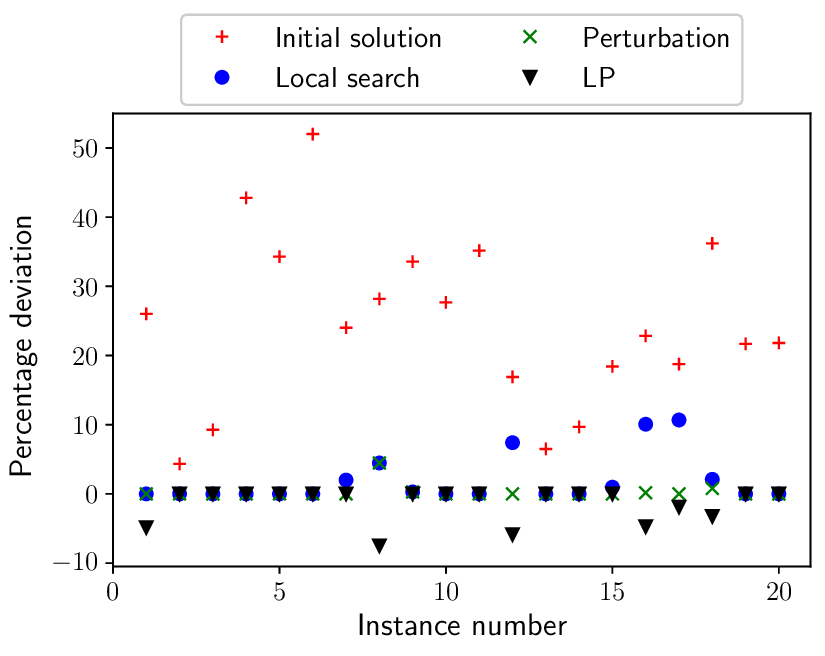}
         \caption{Deviation of heuristic and LP from optimum for $20\%$ target allocation}
         \label{subfig: deviation_heuristic_20_frac_targ_alloc_vehicles_same_depot}
     \end{subfigure}
     \hfill
     \begin{subfigure}[b]{0.46\textwidth}
         \centering
         \includegraphics[width=\textwidth]{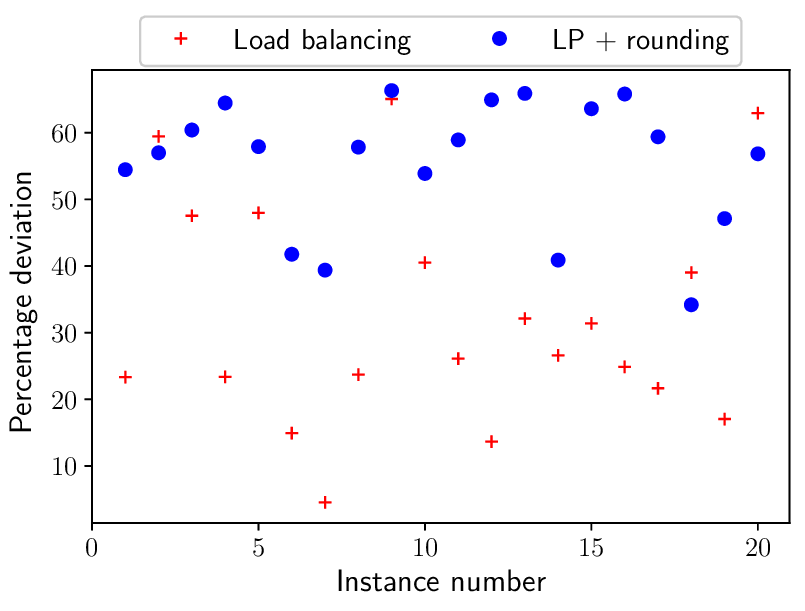}
         \caption{Deviation of two initial feasible solutions from optimum for $0\%$ target allocation}
         \label{subfig: deviation_initial_solution_0_frac_targ_alloc_vehicles_same_depot}
     \end{subfigure}
     \hfill
     \begin{subfigure}[b]{0.46\textwidth}
         \centering
         \includegraphics[width=\textwidth]{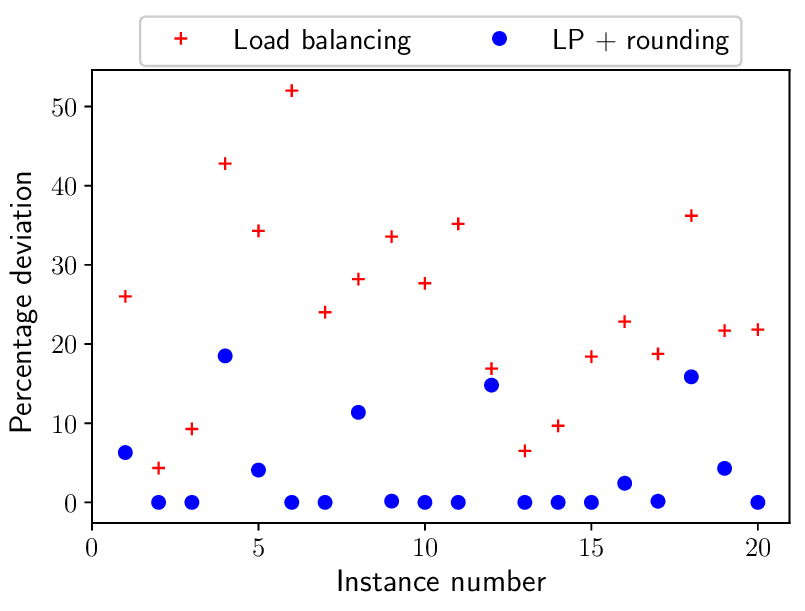}
         \caption{Deviation of two initial feasible solutions from optimum for $20\%$ target allocation}
         \label{subfig: deviation_initial_solution_20_frac_targ_alloc_vehicles_same_depot}
     \end{subfigure}
     \hfill
     \begin{subfigure}[b]{0.46\textwidth}
         \centering
         \includegraphics[width=\textwidth]{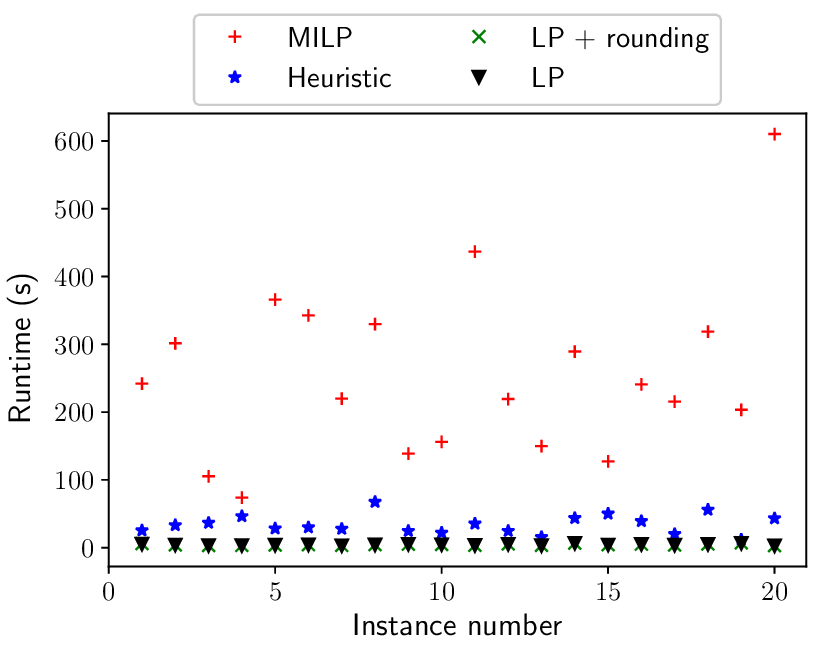}
         \caption{Runtime of branch and cut, heuristic, LP, and LP with rounding for $0\%$ target allocation}
         \label{subfig: runtime_0_frac_targ_alloc_vehicles_same_depot}
     \end{subfigure}
     \hfill
     \begin{subfigure}[b]{0.46\textwidth}
         \centering
         \includegraphics[width=\textwidth]{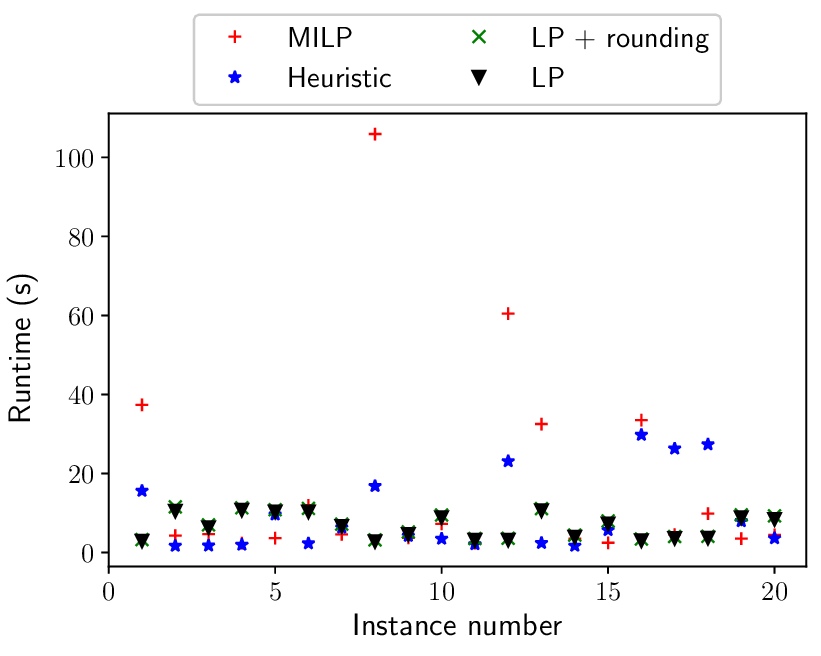}
         \caption{Runtime of branch and cut, heuristic, LP, and LP with rounding for $20\%$ target allocation}
         \label{subfig: runtime_20_frac_targ_alloc_vehicles_same_depot}
     \end{subfigure}
    \caption{Scenario 2: Vehicles $1$ and $2$ have unit speeds and start at same depot, vehicle $3$ has $2$ units speed, thirty targets, twenty instances}
    \label{fig: two_vehicles_start_at_same_depot}
\end{figure}

In this scenario, the speeds of the three vehicles were chosen as $1, 1,$ and $2$ units, respectively. Further, the first and second vehicles were considered to start at the same depot location. The first and third vehicles' depot locations were randomly generated on the $200\times200$ grid. The comparison of the heuristic with the optimal solution obtained using branch and cut for $0\%$ and $20\%$ target assignment and the runtimes are shown in Fig.~\ref{fig: two_vehicles_start_at_same_depot} for twenty instances. Further, the mean percentage deviation from the optimum for the heuristic, the LP, and the initial feasible solutions are summarized in Table~\ref{tab: mean_percentage_deviation_two_vehicles_same_depot}, and the mean running times are summarized in Table~\ref{tab: mean_running_time_scenario_2}. It should be noted that for the cases with $0\%, 10\%,$ and $20\%$ target assignment, the number of instances for which the solution obtained from the heuristic was within $2\%$ of the optimal solution was $11, 13,$ and $19,$ respectively. Other observations regarding the solution quality from the heuristic and LP with rounding, the deviation of the lower bound from the LP, and computation times are similar to the observations for scenario~1.

% \begin{table}[htb!]
%     \centering
%     \begin{tabular}{|c|ccccc|c|}
%     \hline
%     \multirow{3}{*}{\textbf{\begin{tabular}[c]{@{}c@{}}\% target \\ assignment\end{tabular}}} & \multicolumn{5}{c|}{\textbf{Mean percentage deviation from optimum}} & \textbf{No. of instances} \\ \cline{2-6}
%      & \multicolumn{1}{c|}{LP with} & \multicolumn{1}{c|}{\multirow{2}{*}{Load balancing}} & \multicolumn{1}{c|}{\multirow{2}{*}{Local search}} & \multicolumn{1}{c|}{\multirow{2}{*}{Perturbation}} & Max. deviation & \textbf{within $2\%$ of} \\
%      & \multicolumn{1}{c|}{rounding} & \multicolumn{1}{c|}{} & \multicolumn{1}{c|}{} & \multicolumn{1}{c|}{} & of heuristic & \textbf{optimum} \\ \hline
%     0 & \multicolumn{1}{c|}{} & \multicolumn{1}{c|}{} & \multicolumn{1}{c|}{} & \multicolumn{1}{c|}{} &  &  \\ \hline
%     10 & \multicolumn{1}{c|}{} & \multicolumn{1}{c|}{} & \multicolumn{1}{c|}{} & \multicolumn{1}{c|}{} &  &  \\ \hline
%     20 & \multicolumn{1}{c|}{$4.66$} & \multicolumn{1}{c|}{$19.27$} & \multicolumn{1}{c|}{$1.87$} & \multicolumn{1}{c|}{$0.69$} & $8.18$ & $27$ \\ \hline
%     \end{tabular}
%     \caption{Mean percentage deviation of heuristic and initial feasible solution from LP with rounding for $v_1 = v_2 = 1$, $d_1 = d_2$, and $v_3 = 2$, thirty targets, and thirty instances}
% \label{tab: mean_percentage_deviation_two_vehicles_same_depot}
% \end{table}

% \begin{figure}
%      \centering
%      \begin{subfigure}[b]{0.45\textwidth}
%          \centering
%          \includegraphics[width=\textwidth]{Figures/Instance_1_initial_solution_LP_rounding.eps}
%          \caption{$y=x$}
%          \label{fig:y equals x}
%      \end{subfigure}
%      \hfill
%      \begin{subfigure}[b]{0.45\textwidth}
%          \centering
%          \includegraphics[width=\textwidth]{Figures/Instance_1_initial_solution_heuristic.eps}
%          \caption{$y=3\sin x$}
%          \label{fig:three sin x}
%      \end{subfigure}
%      \hfill
%      \begin{subfigure}[b]{0.45\textwidth}
%          \centering
%          \includegraphics[width=\textwidth]{Figures/Instance_1_local_search.eps}
%          \caption{$y=5/x$}
%          \label{fig:five over x}
%      \end{subfigure}
%      \hfill
%      \begin{subfigure}[b]{0.45\textwidth}
%          \centering
%          \includegraphics[width=\textwidth]{Figures/Instance_1_perturbation.eps}
%          \caption{$y=5/x$}
%          \label{fig:five over x}
%      \end{subfigure}
%      \hfill
%      \begin{subfigure}[b]{0.45\textwidth}
%          \centering
%          \includegraphics[width=\textwidth]{Figures/Instance_1_MILP.eps}
%          \caption{$y=5/x$}
%          \label{fig:five over x}
%      \end{subfigure}
%         \caption{Three simple graphs}
%         \label{fig:three graphs}
% \end{figure}

\begin{table}[htb!]
    \centering
    \begin{tabular}{|c|ccccc|c|}
    \hline
    \multirow{3}{*}{\textbf{\begin{tabular}[c]{@{}c@{}}\% target \\ allocation\end{tabular}}} & \multicolumn{5}{c|}{\textbf{Mean percentage deviation from optimum}} & \multirow{3}{*}{\textbf{\begin{tabular}[c]{@{}c@{}}Max. deviation \\ of heuristic from \\ optimum ($\%$)\end{tabular}}} \\ \cline{2-6} 
     & \multicolumn{1}{c|}{LP with} & \multicolumn{1}{c|}{\multirow{2}{*}{Load balancing}} & \multicolumn{1}{c|}{\multirow{2}{*}{Local search}} & \multicolumn{1}{c|}{\multirow{2}{*}{Perturbation}} & \multirow{2}{*}{LP} &  \\
     & \multicolumn{1}{c|}{rounding} & \multicolumn{1}{c|}{} & \multicolumn{1}{c|}{} & \multicolumn{1}{c|}{} & &  \\ \hline
    $0$ & \multicolumn{1}{c|}{$55.55$} & \multicolumn{1}{c|}{$32.29$} & \multicolumn{1}{c|}{$10.43$} & \multicolumn{1}{c|}{$3.65$} & \multicolumn{1}{c|}{$-21.75$} & $12.76$ \\ \hline
    $10$ & \multicolumn{1}{c|}{$27.27$} & \multicolumn{1}{c|}{$27.22$} & \multicolumn{1}{c|}{$6.53$} & \multicolumn{1}{c|}{$2.68$} & \multicolumn{1}{c|}{$-6.72$} & $12.41$ \\ \hline
    $20$ & \multicolumn{1}{c|}{$3.90$} & \multicolumn{1}{c|}{$24.51$} & \multicolumn{1}{c|}{$1.90$} & \multicolumn{1}{c|}{$0.29$} & $-1.41$ & $4.47$ \\ \hline
    \end{tabular}
    \caption{Mean percentage deviation of stages of heuristic, LP, and LP with rounding from optimum for $v_1 = v_2 = 1$, $d_1 = d_2$, and $v_3 = 2$, thirty targets, and twenty instances}
    \label{tab: mean_percentage_deviation_two_vehicles_same_depot}
\end{table}

\begin{table}[htb!]
    \centering
    \begin{tabular}{|c|c|c|c|c|}
    \hline
    \textbf{\% target allocation} & \textbf{Branch and cut} (s) & \textbf{Heuristic} (s) & \textbf{LP} (s) & \textbf{LP with rounding} (s) \\ \hline
    $0$ & $254.38$ & $34.03$ & $4.04$ & $4.52$ \\ \hline
    $10$ & $121.60$ & $24.69$ & $5.83$ & $6.41$ \\ \hline
    $20$ & $17.15$ & $9.69$ & $6.46$ & $7.07$ \\ \hline
    \end{tabular}
    \caption{Mean running time of branch and cut, heuristic, LP, and LP with rounding for $v_1 = v_2 = 1$, $d_1 = d_2$, and $v_3 = 2$, thirty targets, and twenty instances}
    \label{tab: mean_running_time_scenario_2}
\end{table}

\section{Conclusion}

In this paper, a heuristic was proposed for a min-max heterogeneous multi-vehicle multi-depot traveling salesman problem. The vehicles were considered functionally heterogeneous due to different vehicle-target assignments and structurally heterogeneous due to different vehicle speeds. The proposed heuristic extended the MD algorithm, an algorithm that yields good feasible solutions for the homogeneous variant of the considered problem. The proposed heuristic consisted of an initialization stage based on load balancing, a local search stage to improve the maximal tour, and a perturbation stage to break from a local optimum. The heuristic was implemented and benchmarked with the optimal solution obtained by solving a mixed-integer linear program for the problem using branch and cut for instances considering three vehicles and thirty targets. It was observed that the proposed heuristic yielded good feasible solutions, which was about $4\%$ on average from the optimal solution. Future work can include exploring a better feasible solution and utilizing more neighborhoods to improve the quality of the solution obtained.

\section*{Acknowledgments}
This work was funded under FA8650-19-C-1692. This work was cleared for public release under AFRL-2023-3141 and the views or opinions expressed in this work are those of the authors and do not reflect any position or opinion of the United States Government, US Air Force, or Air Force Research Laboratory. 

\bibliography{citations}

\end{document}